\documentclass[a4paper,reqno,10pt]{amsart}
\usepackage{amssymb}

\makeatletter
\renewcommand\eqref[1]{%
	\textcolor{red}{\tagform@{\ref{#1}}}%
}

\usepackage[height=20cm,width=17cm,left=2cm,right=1.5cm]{geometry}
\usepackage{amsfonts}
\usepackage{bigints}
\usepackage{amsmath}
\usepackage{cite}
\usepackage{amssymb}
\usepackage{systeme}
\usepackage{mathrsfs}
\usepackage[colorlinks=true, linktocpage=false]{hyperref}
\usepackage[nameinlink,capitalize]{cleveref}
\numberwithin{equation}{section}
\newtheorem{theorem}{Theorem}[section]

\newtheorem{lm}[theorem]{Lemma}

\newtheorem{rem}{Remark}
\newcommand{\R}{\mathbb{R}}

\newcommand{\N}{\mathbb{N}}

\def\ds{\displaystyle}

\def\dive{\mathrm{div}}
\newenvironment{preuve}{{\noindent {\bf Proof. }}}{\hfill {\rule{2.5mm}{2.5mm}} }

\author[M.~Amara]{Mustapha Amara}
\author[J.~Benameur]{Jamel Benameur}
\address{Department of Mathematics, Faculty of Science of Gab\`es, Research Laboratory Mathematics and Applications LR17ES11; Tunisia}
\email{\sl jamelbenameur@gmail.com}
\email{\sl Mostafa.Amara@fsg.u-gabes.tn}

\title[Existence and Uniqueness of the Solution to the (AQG)
Equations in the Sobolev Space]
{Existence and Uniqueness of the Solution to the Anisotropic Quasi-Geostrophic
	Equations in the Sobolev Space}

\begin{document}
\begin{abstract}
	In this paper, we focus on the two-dimensional surface quasi-geostrophic
	equation with fractional horizontal dissipation and vertical thermal diffusion which represents a general case of the classical surface quasi-geostrophic equation.
	On the one hand, we will show the local existence and uniqueness of the solution in Sobolev space $H^{2-2\alpha}(\R^2)\cap H^{2-2\beta}(\R^2)$, which is the critical space in the classical case.  Furthermore, we will demonstrate that the solution is global even when the initial data is very small. Finally, we will study the asymptotic representation of our global solution in infinity.
\end{abstract}



\maketitle
\tableofcontents


\section{\bf Introduction}
We consider the two-dimensional anisotropic quasi-geostrophic  equation:
	\begin{equation}\label{AQG}\tag*{(AQG)}
		\begin{cases}
			\partial_t\theta+ u_\theta.\nabla\theta +\mu|\partial_1|^{2\alpha}\theta+\nu |\partial_2|^{2\beta}\theta=0,&t>0,x\in\R^2,\\
			u_\theta=\mathcal{R}^\perp\theta &t>0,x\in\R^2,\\
			\theta(x,0)=\theta^0(x),& x\in \R^2,
		\end{cases}
	\end{equation}
\begin{itemize}
	\item[$\bullet$] $\theta=\theta(x,t)$ and $u_\theta=u_\theta(x,t)=(u_\theta^1(x,t),u_\theta^2(x,t))$ are the unknown potential temperature and the velocity field of the fluid, respectively.
	\item[$\bullet$] $\theta^0$ is the given initial potential temperature.
	\item[$\bullet$] $\mu,\nu>0$  and $\alpha,\beta\in (0,1)$. 
	\item[$\bullet$] $\mathcal{R}^\perp=\left(-\mathcal{R}_2,\mathcal{R}_1\right)$, where $\mathcal{R}_1$ and $\mathcal{R}_2$ are the Riesz transforms
	defined by
	$$\mathcal{R}_j \theta=\mathcal{F}^{-1}\left(\xi\mapsto \frac{i\xi_j}{|\xi|}\mathcal{F}(\theta)(\xi)\right),\quad j\in\{1,2\}.$$
\end{itemize}

 The  equation arises from geophysics and has a strong physical background. If $\mu=\nu=0$, then it reduces to the inviscid case, which shares some important features with the 3D Euler equations such as the vortex stretching mechanism. This inviscid case equation is an important model of geophysical fluid dynamics, which describes the evolution of the surface temperature field in the rotating stratified fluid. The first mathematical studies of this equation was carried out in 1994s by Constantin, Majda and Tabak. For more details and mathematical and physical explanations of this model we can consult \cite{CP,CP1,DC,JP}.\\

In this manuscript, we show the local existence and uniqueness of the solution of system \ref{AQG}. More precisely,
for the given initial data $\theta^0$  in $H^{\max\{2-2\alpha,2-2\beta\}}(\R^2):=H^{2-2\alpha}(\R^2)\bigcap H^{2-2\beta}(\R^2) $, there exists a positive
constant $T_0$ determined by $\alpha,\beta$ and $\theta^0$ such that \ref{AQG} possesses a
unique local classical solution on $C([0, T_0],H^{\max\{2-2\alpha,2-2\beta\}}(\R^2))$. Next, we prove the global existence for a small initial data in the same space, that is to say, a positive constant $\varepsilon$ exists such, if $\|\theta^0\|_{H^{\max\{2-2\alpha,2-2\beta\}}}<\varepsilon$, then, the local solution is global one. Finally we prove some optimal decay results of the global solution of \ref{AQG}  with small initial data.\\

Before we state the main result, we recall some known results for the existence theorem of the classical quasi-geostrophic equation:
\begin{equation}\label{SQG}\tag*{(SQG)}
	\begin{cases}
		\partial_t\theta+ u_\theta.\nabla\theta +\kappa(-\Delta)^\gamma\theta=0,&t>0,x\in\R^2,\\
		u_\theta=\mathcal{R}^\perp\theta, &t>0,x\in\R^2,\\
		\theta(x,0)=\theta^0(x),& x\in \R^2.
	\end{cases}
\end{equation}
Easy to see that when $\gamma=\alpha=\beta$ and $\kappa=\mu =\nu$ , we have \ref{AQG}  becomes the classical dissipation $\ref{SQG}$ equation. This last equation was studied in three different cases: subcritical case when $\gamma\in\left(\frac{1}{2},1\right)$, critical case $\gamma=\frac{1}{2}$ and last case when $\gamma\in\left(0,\frac{1}{2}\right)$ called supercritical case.
Constantin and Wu \cite{CW} established the existence of a unique global solution and decay estimates with regard to the $L^2$ norm for the initial data $\theta^0 \in L^2(\R^2)$ in the subcritical case. For the critical case, there are many results for this case, for example, Conctantin and Vicol \cite{CV} showed the global in time existence of the smooth solution $\theta\in \mathcal{S}(\R^2)$.  See also the works \cite{CZ,ZF,KNV} where same type of results have
been obtained.\\

 However, the supercritical cases, for an initial data, whether it remains globally regular or not, is an interesting open problem. Although the global well-posedness is still open for the this case for \ref{SQG} equation. Some interesting result, showed by Miura \cite{MH} who proved the unique local existence of the solution  in the critical space $H^{2-2\gamma}(\R^2)$.  See also the works in \cite{JN} that have produced similar results.\\

To our knowledge, the first time the system of equations in \ref{AQG} has been studied, is by Zhuan in \cite{YZ}, who established the global regularity when the dissipation powers are restricted to
a suitable range. More specifically
\begin{theorem}[see \cite{YZ}]\label{theorem1.1}
	Let $\theta^0\in H^s(\R^2)$ with $s\geq 2$ and $\alpha, \beta>0$. Then there exists a positive $T(\|\theta^0\|_{H^s})>0$ such that for the system $(AQG)$
	admits a unique solution $$\theta\in C([0,T],H^s(\R^2)),\quad|\partial_{1}|^\alpha\theta,|\partial_{2}|^\beta\theta\in L^2([0,T],H^s(\R^2)).$$
	Moreover, if $\alpha,\beta\in (0,1)$ satisfy 
	\begin{equation}\label{1.1}
		\beta>\begin{cases}
			\frac{1}{2\alpha+1},&0<\alpha\leq \frac{1}{2}\\
			\\
			\frac{1-\alpha}{2\alpha},& \frac{1}{2}<\alpha<1,
		\end{cases}
	\end{equation}
 then, the system $(AQG)$ admits a unique global solution $\theta$ such that for any $T>0$ we have
	$$
	\theta\in C([0,T],H^s(\R^2)),\ |\partial_{1}|^{\alpha}\theta,|\partial_{2}|^{\beta}\theta\in L^2([0,T],H^s(\R^2)).
	$$
\end{theorem}
Motivated by these previous studies, we first show the local solution existence to \ref{AQG} for the initial data $\theta^0\in H^{s}(\R^2)$, with $s=\max\{2-2\alpha,2-2\beta\}$. Our main result is as follows:
\begin{theorem}\label{theorem1.2} Let $\alpha,\beta\in (0,1)$ such that
	$$\min\{\alpha,\beta\}<\frac{1}{2}.$$
	We assume $s=\max\{2-2\alpha,2-2\beta\}$. Let $\theta^0\in H^{s}(\R^2)$, then, there exists
	a positive time $T_0$  such that \ref{AQG}  admits a unique solution
	$$\theta\in C([0,T_0];H^{s}(\R^2)),\ |\partial_1|^{\alpha}\theta,|\partial_2|^{\beta}\theta\in L^2([0,T_0];H^{s}(\R^2)).$$
	Moreover, there exists a constant $c>0$ such that, if $$\|\theta^0\|_{\dot{H}^{s}}<c,$$ then
	$$\theta\in C(\R^+;H^{s}(\R^2)),\ |\partial_1|^{\alpha}\theta,|\partial_2|^{\beta}\theta\in L^2(\R^+;H^{s}(\R^2)).$$
	In addition the following estimate holds
	\begin{equation}\label{1.2}
		\|\theta(t)\|_{H^{s}}^2+\int_{0}^{t}\||\partial_{1}|^\alpha\theta(\tau)\|_{H^{s}}^2 d\tau+\int_{0}^{t}\||\partial_{2}|^\beta\theta(\tau)\|_{H^{s}}^2 d\tau\leq \|\theta^0\|_{H^{s}}^2.
	\end{equation}
\end{theorem}
\begin{rem}
\begin{itemize}
	\item[]
	\item[{\it (i)}] 	We assume that $\mu=\nu=1$ to simplify the calculus and some procedures in the proofs of our results.
	\item[{\it (ii)}] We can select without loss the generality $\alpha=\min\{\alpha,\beta\},$ the proof will be in the same way if we have the opposite.
	\item[{\it (iii)}] Very recently in \cite{YZ1}, Zhuan a proved the global well-posedness existence of \ref{AQG} in the Sobolev space
	$ H^{\max{s,2-\frac{4\alpha\beta}{\alpha+\beta}}}:=H^{s}(\R^2)\cap \dot{H}^{2-\frac{4\alpha\beta}{\alpha+\beta}}(\R^2)$, with $s\geq 0$, when $$\|\theta^0\|_{\dot{H}^{2-\frac{4\alpha\beta}{\alpha+\beta}}}<\varepsilon,$$
	where $\varepsilon>0$ is very small. This result not true only if $(\alpha,\beta)\in \mathbb{F}$, where 
	$$\mathbb{F}:=(0,1]^{2}\setminus\left\{(\alpha,\beta)\in (0,1]^{2} \mbox{ satisfy \eqref{1.1}} \right\}.$$
	In the same paper, the author showed the decay estimate for the solution with additionally that $$\theta^0\in H^{2-\frac{4\alpha\beta}{\alpha+\beta}}(\R^2)\cap L^p(\R^2),\quad p\in [1,2).$$
\end{itemize}
\end{rem}

	In the next theorem for the decays as time goes to infinity of the solution:
	\begin{theorem}\label{theorem1.3}
		Let $\theta\in C(\R^+,H^{s}(\R^2))$, $s=\max\{2-2\alpha,2-2\beta\}$, a global solution of \ref{AQG} equation such that $$\|\theta^0\|_{\dot{H}^{s}}<c,$$ then
		\begin{equation*}
			\lim_{t\rightarrow +\infty}\|\theta(t)\|_{H^{s}}=0.	
		\end{equation*}
	\end{theorem}

The following is a breakdown of the paper's structure. In Section 2, we review some of the notations and estimations that will be used in the next parts. In Section 3, we prove the existence theorem for any initial data $\theta^0\in H^{2-2\alpha}(\R^2)$. Finally, in Section 4, we prove \cref{theorem1.3}.\\

Throughout this paper, the constant, which may differ in each line, is denoted by $C$ throughout this work. $C = C(a_1,...,a_n)$ specifies that $C$ is entirely dependent on $a_1,...,a_n$.

 	\section{\bf Notations and Preliminary Results}
 In this short section, we collect some notations and definitions that will be used later, and we	give some technical lemmas.
 \subsection{\bf Notations}
 Let $\mathcal{S}(R^2)$ denote the space of Schwartz class functions defined on $\R^2$ and $\mathcal{S}'(R^2)$ denote the space of tempered distributions. For $f\in\mathcal{S}(R^2)$ , we denote by $\widehat{f}$ or $\mathcal{F}(f)$, the Fourier transform of $f$, defined as
 	\begin{equation}
 	\mathcal{F}(f)(\xi)=\widehat{f}(\xi)=\int_{\R^2}e^{-ix.\xi}f(x)dx,\quad \forall\xi\in \R^2,
 \end{equation}
and the inverse Fourier transform of $f$ defined by
\begin{equation}
	\mathcal{F}^{-1}(f)(x)=(2\pi)^{-2}\int_{\R^2}e^{i\xi.x}f(\xi)d\xi,\quad \forall x\in \R^2.
\end{equation}
Recall that $\mathcal{F}$ is an isometry on $L^2$ and satisfies 
\begin{equation}
	(f,g)_{L^2}=(\widehat{f},\widehat{g})_{L^2}.
\end{equation}
$\bullet$ The convolution product of a suitable pair of function $f=(f_1,f_2)$ and $g=(g_1,g_2)$ on $\R^2$ 	is given by
$$
f\ast g(x)=\int_{\R^2} f(x-y)g(y)dy.
$$
Also, we set
$$f\otimes g:= (g_1f,g_2f)\quad\mbox{ and }\quad\dive (f\otimes g):= (\dive(g_1f),\dive(g_2f)).$$
$\bullet$ The fractional operators are defined through the Fourier transform, namely, for any $s\in \R$
\begin{equation}
	|\partial_1|^{s}f=\mathcal{F}^{-1}\left(\xi\mapsto|\xi_1|^{s}\widehat{f}(\xi)\right),\; |\partial_2|^{s}f=\mathcal{F}^{-1}\left(\xi\mapsto|\xi_2|^{s}\widehat{f}(\xi)\right)\mbox{ and }
|\nabla|^{s}f=\mathcal{F}^{-1}\left(\xi\mapsto|\xi|^{s}\widehat{f}(\xi)\right).
\end{equation}
$\bullet$, For $s\in \R$, let $H^s(\R^2)$ denote the Sobolev non-homogeneous space of order $s$, which is a Hilbert space with the
inner product
$$\left(f,g\right)_{H^s}=\left((1+|\nabla|^{2})^{\frac{s}{2}}f,(1+|\nabla|^{2})^{\frac{s}{2}}g\right)_{L^2}=\int_{\R^2}(1+|\xi|^2)^s\widehat{f}(\xi)\overline{\widehat{g}(\xi)}d\xi$$
and the norm $\|f\|_{H^s}=\ds\sqrt{\left(f,f\right)_{H^s}}$. We also denote $\dot{H}^s(\R^2)$ the Sobolev homogeneous space of order $s$ with the
inner product
$$\left(f,g\right)_{\dot{H}^s}=\left(|\nabla|^{s}f,|\nabla|^{s}g\right)_{L^2}=\int_{\R^2}|\xi|^{2s}\widehat{f}(\xi)\overline{\widehat{g}(\xi)}d\xi$$
and the norm $\|f\|_{\dot{H}^s}=\ds\sqrt{\left(f,f\right)_{\dot{H}^s}}$. If $s > 0$, then it is well known that the norm equivalence
$$\|f\|_{H^s}\sim \|f\|_{L^2}+\|f\|_{\dot{H}^s}.$$
$\bullet$ If $(X, \|.\|_{X})$, be a Banach space, $p\in [0,+\infty]$ and $T > 0$. We define $$L_T^p(X):=L^p([0,T],X)$$ the space of all
measurable functions $t\in[0,T]\mapsto f(t)\in X $ such that 
$$\left(t\mapsto \|f(t)\|_{X}\right)\in L^p([0,T]).$$
$\bullet$ We present some features for the mollifier in $\R^2$ of the Friedrichs type, as specified by
$$J_N(f):=\mathcal{F}^{-1}\left(\xi\mapsto \mathcal{X}_{B(0,N)}(\xi)\widehat{f}(\xi)\right),$$
where $N\in \N$, $f\in L^2(\R^2)$, $B(0,N)=\{\xi\in \R^2;\ |\xi|<N\}$ and 
\begin{align*}
	\begin{array}{rccl}
		\mathcal{X}_{B(0,N)}:&\R^2&\rightarrow& \{0,1\}\\
		&\xi&\mapsto &\mathcal{X}_{B(0,N)}(\xi)=\begin{cases}
			1&\mbox{ if }|\xi|<1,\\
			0&\mbox{ else. }
		\end{cases}
	\end{array}
\end{align*}
\subsection{\bf Preliminary Results}
Next, we introduce some Lemma that will be used in the proof of our results.
 \begin{lm}[see \cite{HB}]\label{Lemma2.1}
 	Let $H$ be Hilbert space and $(x_n)$ be a bounded sequence
 	of elements in $H$ such that
 	\begin{align*}
 		&x_n\rightharpoonup x\mbox{ in }H \quad\mbox{ and }\quad \lim_{n\rightarrow+\infty}\|x_n\|\leq \|x\|.
 	\end{align*}
 Therefore $\ds\lim_{n\rightarrow+\infty}\|x_n-x\|=0$.
 \end{lm}
 \begin{lm}[see \cite{BH}]\label{Lemma2.2}
 	Let $p \in [2, +\infty)$ and $\sigma\in [0,1)$ such that
 	$$\frac{1}{p}+\frac{\sigma}{2}=\frac{1}{2}.$$
 	Then, there is a constant $C > 0$ such	 that
 	$$\|f\|_{L^p(\R^d)}\leq C\||\nabla|^{\sigma}f\|_{L^2(\R^d)}.$$
 \end{lm}
 \begin{lm}[see \cite{BH}]\label{Lemma2.3}
 	Let $s_1$, $s_2$ be two real numbers such that $s_1<1$ and $s_1+s_2>0$. Then, there exists a positive constant $C=C(s_1,s_2)$ such that for all $f,g\in \dot{H}^{s_1}(\R^2)\bigcap \dot{H}^{s_2}(\R^2)$; 
 	\begin{equation}
 		\|fg\|_{\dot{H}^{s_1+s_2-1}}\leq C \left(\|f\|_{\dot{H}^{s_1}}\|g\|_{\dot{H}^{s_2}}+\|f\|_{\dot{H}^{s_2}}\|g\|_{\dot{H}^{s_1}}\right).
 	\end{equation}
 	Moreover, in addition $s_2<1$, there exists a positive constant $C=C(s_1,s_2)$ such that, for all $f\in \dot{H}^{s_1}(\R^2)$ and $g\in \dot{H}^{s_2}(\R^2)$; 
 	\begin{equation}
 		\|fg\|_{\dot{H}^{s_1+s_2-1}}\leq C\|f\|_{\dot{H}^{s_1}}\|g\|_{\dot{H}^{s_2}}
 	\end{equation}
 \end{lm}
We recall Sobolev's interpolation theorem:
\begin{lm}\label{Lemma2.6}
	For $s_1,s_2\in\R$ and $t\in[0,1]$; the interpolation inequalities, respectively, in the homogeneous and non-homogeneous Sobolev spaces
	\begin{align}
		&\|f\|_{H^{ts_1+(1-t)s_2}}\leq \|f\|_{H^{s_1}}^t\|f\|_{H^{s_2}}^{1-t},\\
		&\|f\|_{\dot{H}^{ts_1+(1-t)s_2}}\leq \|f\|_{\dot{H}^{s_1}}^t\|f\|_{\dot{H}^{s_2}}^{1-t}.
	\end{align}
\end{lm}
	\begin{lm}\label{Lemma2.7}
	Let $s,s'\in \R$ and $\alpha,\beta\in (0,1)$, such that $s'<s+\alpha$ and $\alpha\leq \beta$ then for any $f\in \mathcal{S}(\R^2)$
	\begin{equation*}
		\||\nabla|^\alpha f\|_{\dot{H}^{s}}\leq \| f\|_{\dot{H}^{s'}} + \||\partial_{1}|^\alpha f\|_{\dot{H}^s}+ \||\partial_{2}|^\beta f\|_{\dot{H}^s}.
	\end{equation*}
\end{lm}
\begin{preuve}
	For $f\in \mathcal{S}(\R^2)$, we have 
	\begin{align*}
		\||\nabla|^\alpha f\|_{\dot{H}^s}^2&=\int_{\R^2}(|\xi_1|^2+|\xi_2|^2)^{\alpha}\left|\mathcal{F}(|\nabla|^{s}f)(\xi)\right|^2d\xi\\
		&\leq \int_{\R^2}|\xi_1|^{2\alpha}\left|\mathcal{F}(|\nabla|^{s}f)(\xi)\right|^2d\xi+\int_{\R^2}|\xi_2|^{2\alpha}\left|\mathcal{F}(|\nabla|^{s}f)(\xi)\right|^2d\xi=\||\partial_{1}|^\alpha f\|_{\dot{H}^s}^2+\||\partial_{2}|^\alpha f\|_{\dot{H}^s}^2.
	\end{align*}
	But $s'-s<\alpha\leq \beta$ then their exist $z\in (0,1]$ such that $\alpha=z\times\beta+(1-z)\times(s'-s)$ and by interpolation inequality we get
	\begin{align*}
		\||\partial_{2}|^\alpha f\|_{\dot{H}^s}&\leq \| f\|_{\dot{H}^{s'}}^{1-z}\||\partial_{2}|^\beta f\|_{\dot{H}^s}^{z}\\
		&\leq (1-z)\| f\|_{\dot{H}^{s'}}+z\||\partial_{2}|^\beta f\|_{\dot{H}^s} \leq \| f\|_{\dot{H}^{s'}}+\||\partial_{2}|^\beta f\|_{\dot{H}^s}.
	\end{align*}
\end{preuve}
 \begin{lm}[see \cite{JN}]\label{Lemma2.4}
	For any $p\in (1,+\infty)$, there is a constant $C(p)>0$ such that
	\begin{equation}
		\|\mathcal{R}^\perp\theta\|_{L^p}\leq C(p) \|\theta\|_{L^p}.
	\end{equation}
\end{lm}
\par We recall the following important commutator and product estimates:
\begin{lm}\label{Lemma2.5}
	For $s>1$, if $f,g\in \mathcal{S}(\R^2)$ then for any $\alpha\in(0,1)$
	\begin{align}
		\label{2.8}	\||\nabla|^s(f g)-f|\nabla|^s g\|_{L^2}\leq s2^{s}C(\alpha)\left(\||\nabla|^{s+\alpha} f\|_{L^{2}}\||\nabla|^{1-\alpha}g\|_{L^{2}}+\||\nabla|^{s-1+\alpha} g\|_{L^{2}}\||\nabla|^{2-\alpha}f\|_{L^{2}}\right).
	\end{align}
\end{lm}
\begin{preuve} We have
	\begin{align*}
		\||\nabla|^s(f g)-f|\nabla|^s g\|_{L^2}^2&=\int_{\R^2}\left|\mathcal{F}(|\nabla|^s(f g)-f|\nabla|^s g)(\xi)\right|^2 d\xi\\
		&\leq \int_{\R^2}\left(\int_{\R^2}\left||\xi|^s-|\eta|^s\right||\widehat{f}|(\xi-\eta)|\widehat{g}|(\eta)d\eta\right)^2 d\xi.
	\end{align*}
	By using the elementary inequality
	$$\left||\xi|^s-|\eta|^s\right|\leq s2^{s-1}\left(|\xi-\eta|^s+|\eta|^{s-1}|\xi-\eta|\right)$$
	we get
	\begin{align*}
		\||\nabla|^s(f g)-f|\nabla|^s g\|_{L^2}^2&\leq s^2 2^{2(s-1)} \int_{\R^2}\left(\int_{\R^2}|\xi-\eta|^s|\widehat{f}(\xi-\eta)||\widehat{g}(\eta)|+|\xi-\eta||\widehat{f}(\xi-\eta)||\eta|^{s-1}|\widehat{g}(\eta)|d\eta\right)^2 d\xi\\
		&\leq s^22^{2s}\|f_1g_1\|^2_{L^2}+s^22^{2s}\|f_2g_2\|^2_{L^2},
	\end{align*}
	where
	\begin{align*}
		&\mathcal{F}(f_1)(\xi)=|\xi|^s|\mathcal{F}(f)(\xi)|&& \mathcal{F}(g_1)(\xi)=|\mathcal{F}(g)(\xi)|\\
		&\mathcal{F}(f_2)(\xi)=|\xi||\mathcal{F}(f)(\xi)|&& \mathcal{F}(g_2)(\xi)=|\xi|^{s-1}|\mathcal{F}(g)(\xi)|\\
	\end{align*}
	Using Hölder inequality with $$\frac{\alpha}{2} +\frac{1-\alpha}{2}=\frac{1}{2}$$
	and by \cref{Lemma2.2}, we get
	\begin{align*}
		\||\nabla|^s(f g)-f|\nabla|^s g\|_{L^2}&\leq s2^s\left( \|f_1\|_{L^{\frac{2}{1-\alpha}}}\|g_1\|_{L^{\frac{2}{\alpha}}}+ \|f_2\|_{L^{\frac{2}{\alpha}}}\|g_1\|_{L^{\frac{2}{1-\alpha}}}\right)\\
		&\leq s2^sC(\alpha) \left( \||\nabla|^{s+\alpha}f\|_{L^{2}}\||\nabla|^{1-\alpha}g\|_{L^{2}}+\||\nabla|^{2-\alpha}f\|_{L^{2}}\||\nabla|^{s-1+\alpha}g\|_{L^{2}}\right),
	\end{align*}
	which finished the proof.
\end{preuve}
\begin{lm}\label{Lemma2.8}
	Let $\alpha,\beta\in (0,1)$ and we pose $A(\xi)=|\xi_1|^{2\alpha}+|\xi_2|^{2\beta}$, then
	\begin{equation}
	|\xi|\leq C(\alpha,\beta)\left(A(\xi)^{\frac{1}{2\alpha}}+A(\xi)^{\frac{1}{2\beta}}\right), \quad\forall \xi\in \R^2.
	\end{equation}
\end{lm}
\begin{preuve}
	Let $\alpha,\beta\in (0,1)$,
	\begin{itemize}
		\item [$\ast$] If $|\xi_1|\leq |\xi_2|$ then we have
		\begin{align*}
			&|\xi|^{2\beta}\leq (|\xi_1|^{2\beta}+|\xi_2|^{2\beta})\leq 2  |\xi_2|^{2\beta}\leq 2 A(\xi).
		\end{align*}
	Therefore
	\begin{align}\label{2.12}
		&|\xi|\leq 2^{\frac{1}{2\beta}}A(\xi)^{\frac{1}{2\beta}}\leq 2^{\frac{1}{2\beta}}\left(A(\xi)^{\frac{1}{2\beta}}+A(\xi)^{\frac{1}{2\alpha}}\right).
	\end{align}
		\item [$\ast$]  If $|\xi_2|\leq |\xi_1|$ then we have
		\begin{align*}
			&|\xi|^{2\alpha}\leq 2(|\xi_1|^{2\alpha}+|\xi_2|^{2\alpha})\leq 2 |\xi_1|^{2\alpha}\leq 2 A(\xi).
		\end{align*}
	Therefore
	\begin{align}\label{2.13}
		&|\xi|\leq 2^{\frac{1}{2\alpha}}A(\xi)^{\frac{1}{2\alpha}}\leq 2^{\frac{1}{2\alpha}}\left(A(\xi)^{\frac{1}{2\alpha}}+A(\xi)^{\frac{1}{2\beta}}\right).
	\end{align}
Collecting \eqref{2.12} and \eqref{2.13}, we get the result.
	\end{itemize}
\end{preuve}
 	\section{\bf Proof of \cref{theorem1.2}}
 	In this section we prove \cref{theorem1.2}. In the first part, we show the local existence and the
 	uniqueness of the solution. The second part is devoted to prove that this solution is continues. The global existence of solution will showed in the last part.
 	\subsection{\bf Local existence and uniqueness:} We start by proving the existence of local solution of \ref{AQG} equation in the following space 
 	$$\mathcal{X}_{T_0}=\left\{f\in L^\infty_{T_0}(H^{2-2\alpha}(\R^2))\bigcap C_{T_0}(L^2(\R^2)):\ |\partial_{1}|^\alpha f,|\partial_{2}|^\beta f\in L^2_{T_0}(H^{2-2\alpha}(\R^2))\right\}.$$
 	
 \vskip0.5cm\noindent{\bf $\bullet$ \underline{\large Existence:}}	The classical study of the existence of the solution necessitates a very small initial date, therefore we will write the initial condition as a sum of higher and lower frequencies.\\
 	
 	 For that, for any $\varepsilon>0$, their exist $N\in \N$ such that $\theta^0=a^0+b^0$, where
 	\begin{equation}
 		\begin{cases}
 			\ast\	a^0:=\mathcal{J}_{N}(\theta^0)\in\ds \bigcap_{r\in \R} H^{r}(\R^2) ,\\ 
 			\ast\	b^0:= \theta^0-a^0,\\ 
 			\ast\	\|b^0\|_{\dot{H}^{2-2\alpha}}<\varepsilon.
 		\end{cases}
 	\end{equation}
According to \cref{theorem1.1}, we get that the following system
 \begin{equation*}
 \begin{cases}
 		\partial_t a+|\partial_{1}|^{2\alpha}a+|\partial_{2}|^{2\beta}a+u_a.\nabla a=0,\\
 		a(0)=a^0\in H^5(\R^2),
 	\end{cases}
 \end{equation*}
admit a unique solution $a\in C_{T^0}(H^5(\R^2))\subset C_{T^0}(H^{2-2\alpha}(\R^2))$, where $T^0=\dfrac{C_N}{\|\theta^0\|_{H^{2-2\alpha}}^4}>0$, moreover,
\begin{align*}
	\|a\|_{L^\infty_{T^0}(H^{5})}\leq M_N.
\end{align*}
To show the existence of solution of system \ref{AQG}, it suffices to show that the following system .	
 	\begin{equation*}\label{(AQG)'}\tag*{(AQG)'}
 			\begin{cases}
 				\partial_t b+|\partial_{1}|^{2\alpha}b+|\partial_{2}|^{2\beta}b+u_a.\nabla b+u_b.\nabla a+u_b.\nabla b=0,&(x,t)\in\R^2\times[0,T^0],\\
 				b(x,0)=b^0(x),&x\in\R^2.
 			\end{cases}
 		\end{equation*}
 	admits a solution in $\mathcal{X}_{T_0}$, where $T_0>0$ is a variable that will be decided later.\\
 	
 	In order to do this, we may for instance
 	make use of the Friedrichs method. The first step is to consider the following approximate system
 	of \ref{(AQG)'},	
 	\begin{equation*}\label{AQGn}\tag*{$(AQG)'_n$}
 	\begin{cases}
 			\partial_tb+|\partial_1|^{2\alpha}\mathcal{J}_nb+|\partial_2|^{2\beta}\mathcal{J}_n b+\mathcal{J}_n\left(u_a.\nabla\mathcal{J}_nb\right)+\mathcal{J}_n\left(\mathcal{J}_nu_b.\nabla a\right)+\mathcal{J}_n\left(\mathcal{J}_nu_b.\nabla\mathcal{J}_nb\right)=0\\
 			b(x,0)=\mathcal{J}_nb^0(x).
 		\end{cases}
 	\end{equation*}
 Using the Cauchy-Lipschitz theorem, for any fixed $n\in \N$, there exist $T_n = T(n, \|b^0\|_{H^{2-2\alpha}})\in(0,T^0)$, such that, the system \ref{AQGn} admit a unique local solution $b_n$ on $[0, T_n]$. Moreover, $\mathcal{J}_nb_n$ is also a solution of \ref{AQGn} with the same initial data, the fact that $\mathcal{J}_n^2=\mathcal{J}_n$. According to the uniqueness, we have
 	$$\mathcal{J}_nb_n=b_n \mbox{ and } \mathcal{J}_nu_{b_n}=u_{b_n}.$$
 	Consequently, for any $n\in \N$, $b_n$ is solution of the following equation
 	\begin{equation}\label{3.2}
 			\partial_tb_n+|\partial_1|^{2\alpha}b_n+|\partial_2|^{2\beta} b_n+\mathcal{J}_n\left(u_a.\nabla b_n\right)+\mathcal{J}_n\left(u_{b_n}.\nabla a\right)+\mathcal{J}_n\left(u_{b_n}.\nabla b_n\right)=0
 	\end{equation}
 	Taking the inner product of \eqref{3.2} with $b_n$ we get
 	\begin{align*}
 		\frac{1}{2}\frac{d}{dt}\|b_n\|_{L^2}^2+\||\partial_{1}|^\alpha b_n\|_{L^2}^2+\||\partial_{2}|^\beta b_n\|_{L^2}^2&\leq \left|\left(\mathcal{J}_n\left(u_{b_n}.\nabla a\right),b_n\right)_{L^2}\right|\\
 		&\leq \|u_{b_n}.\nabla a\|_{L^2}\|b_n\|_{L^2}\\
 		&\leq \| u_{b_n}\|_{L^2}\|\nabla a\|_{L^\infty}\|b_n\|_{L^2}\\
 		&\leq M_N \|b_n\|_{L^2}^2
 	 	\end{align*}
 	We integer in $[0,t]$, $t\in  (0,T_n]$, we get
 	\begin{equation}\label{3.3}
 		\|b_n(t)\|_{L^2}^2+2\int_{0}^t\||\partial_{1}|^\alpha b_n(\tau)\|_{L^2}^2d\tau+2\int_{0}^t\||\partial_{2}|^\beta b_n(\tau)\|_{L^2}^2d\tau\leq \|b^0\|_{L^2}^2+2M_N \int_{0}^t \|b_n(\tau)\|_{L^2}^2d\tau.
 	\end{equation}
 	Gronwall's Lemma implies that
 	\begin{align*}
 		\|b_n(t)\|_{L^2}^2+2\int_{0}^t\||\partial_{1}|^\alpha b_n(\tau)\|_{L^2}^2d\tau+2\int_{0}^t\||\partial_{2}|^\beta b_n(\tau)\|_{L^2}^2d\tau&\leq \|b^0\|_{L^2}^2e^{2M_Nt}\leq \|b^0\|_{L^2}^2e^{2M_NT^0}.
 	\end{align*}
 	Hence, for any $n\in \N$, we have
 	$$b_n\in C^1_{T^0}(H^\sigma(\R^2)),\ |\partial_{1}|^\alpha b_n,|\partial_{2}|^\beta b_n\in L^2_{T^0}(H^\sigma(\R^2)),\quad \forall\sigma\in \R.$$
 	Now, we obtain by applying $|\nabla|^{{2-2\alpha}}$ to \eqref{3.2}
 	and taking $L^2$ inner product
 	with $|\nabla|^{{2-2\alpha}}b_n$ that		
 	\begin{align}\label{3.4}
 		\nonumber	 \frac{1}{2}\frac{d}{dt}\|b_n\|^2_{\dot{H}^{2-2\alpha}}+\||\partial_1|^\alpha b_n\|^2_{\dot{H}^{2-2\alpha}}+\||\partial_2|^\beta b_n\|^2_{\dot{H}^{2-2\alpha}} &=-\left(|\nabla|^{2-2\alpha}(u_a.\nabla b_n),|\nabla|^sb_n\right)_{L^2}\\
 		&\hskip 0.7cm- \left(|\nabla|^{2-2\alpha}(u_{b_n}.\nabla a),|\nabla|^sb_n\right)_{L^2}\\
 		\nonumber	&\hskip 1.4cm- \left(|\nabla|^{2-2\alpha}(u_{b_n}.\nabla b_n),|\nabla|^sb_n\right)_{L^2}\\
 		\nonumber	&	=\mathcal{H}_1+\mathcal{H}_2+\mathcal{H}_3,
 	\end{align}
 where 
 \begin{align*}
 	&\mathcal{H}_1=-\left(|\nabla|^{2-2\alpha}(u_a.\nabla b_n),|\nabla|^sb_n\right)_{L^2}\\
 	&\mathcal{H}_2=-\left(|\nabla|^{2-2\alpha}(u_{b_n}.\nabla a),|\nabla|^sb_n\right)_{L^2}\\
 	&\mathcal{H}_3=- \left(|\nabla|^{2-2\alpha}(u_{b_n}.\nabla b_n),|\nabla|^sb_n\right)_{L^2}.
 \end{align*}
 	In what follows, we shall estimate the terms at the right hand side of \eqref{3.4} one by one.\\
 	To estimate the first term, we use  $\dive u_a=0$ and the first inequality in  \cref{Lemma2.5}  to	conclude
 	\begin{align*}
 		\left|\mathcal{H}_1\right|&=	\left|\left(|\nabla|^{2-2\alpha}(u_a.\nabla b_n),|\nabla|^{2-2\alpha} b_n\right)_{L^2}\right|\\
 		&=\left|\left(|\nabla|^{2-2\alpha}(u_a.\nabla b_n)-u_a.\nabla|\nabla|^{2-2\alpha} b_n,|\nabla|^{2-2\alpha} b_n\right)_{L^2}\right|\\
 		&\leq \||\nabla|^{2-2\alpha}(u_a.\nabla b_n)-u_a.\nabla|\nabla|^{2-2\alpha} b_n\|_{L^2}\||\nabla|^{2-2\alpha}b_n\|_{L^2}\\
 		&\leq C\Big(\||\nabla|^{{2-2\alpha}+\alpha}u_a\|_{L^2}\||\nabla|^{1-\alpha}\nabla b_n\|_{L^2}+\||\nabla|^{{2-2\alpha}-1+\alpha}\nabla b_n\|_{L^2}\||\nabla|^{2-\alpha}u_a\|_{L^2}\Big)\||b_n\|_{\dot{H}^{2-2\alpha}}\\
 		&\leq C_N\||\nabla|^{\alpha} b_n\|_{\dot{H}^{2-2\alpha}}\|b_n\|_{\dot{H}^{2-2\alpha}}.
 	\end{align*}
 	The fact that $\alpha\leq \beta$, then by \cref{Lemma2.7}, we have
 	\begin{align*}
 		\||\nabla|^{\alpha} b_n\|_{\dot{H}^{2-2\alpha}}&\leq \| b_n\|_{L^2}+\||\partial_{1}|^{\alpha} b_n\|_{\dot{H}^{2-2\alpha}}+\||\partial_{2}|^{\beta} b_n\|_{\dot{H}^{2-2\alpha}}\\
 		&\leq \|b^0\|_{L^2}e^{M_NT^0}+\||\partial_{1}|^{\alpha} b_n\|_{\dot{H}^{2-2\alpha}}+\||\partial_{2}|^{\beta} b_n\|_{\dot{H}^{2-2\alpha}},
 	\end{align*}
 	which implies
 	\begin{align}
 		\nonumber	\left|\mathcal{H}_1\right|&\leq C_N\left(1+\||\partial_{1}|^{\alpha} b_n\|_{\dot{H}^{2-2\alpha}}+\||\partial_{2}|^{\beta} b_n\|_{\dot{H}^{2-2\alpha}}\right)\|b_n\|_{\dot{H}^{2-2\alpha}}\\
 		\label{3.5}	&\leq \frac{1}{2} \||\partial_{1}|^{\alpha} b_n\|_{\dot{H}^{2-2\alpha}}^2+\frac{1}{2} \||\partial_{2}|^{\beta} b_n\|_{\dot{H}^{2-2\alpha}}^2+C_N(1+\|b_n\|_{\dot{H}^{2-2\alpha}}^2).
 	\end{align}
 	Move now to estimate the term $\mathcal{H}_2$, the fact that $2-2\alpha>1$ we have 
 	\begin{align}
 		\nonumber	\left|\mathcal{H}_2\right|	&\leq \|u_{b_n}.\nabla a\|_{H^{2-2\alpha}}\|b_n\|_{H^{2-2\alpha}}\\
 		\nonumber	&\leq \|u_{b_n}\|_{H^{2-2\alpha}}\|a\|_{H^{3-2\alpha}}\|b_n\|_{H^{2-2\alpha}}\\
 	\label{3.6}	&\leq C_N(1+\|b_n\|_{\dot{H}^{2-2\alpha}}^2).
 	\end{align}
 	Finally, following the estimate of $\mathcal{H}_1$, with $\dive u_{b_n}=0$, we get
 	\begin{align}
 		\nonumber	\left|\mathcal{H}_3\right|&=\left|\left(|\nabla|^{2-2\alpha}(u_{b_n}.\nabla b_n)-u_{b_n}.\nabla|\nabla|^{2-2\alpha} b_n,|\nabla|^{2-2\alpha}b_n\right)_{L^2}\right|\\
 		\nonumber	&\leq   C\||\nabla|^{{\alpha}}b_n\|_{\dot{H}^{2-2\alpha}}^2\|b_n\|_{\dot{H}^{2-2\alpha}}\\
 	\nonumber	&\leq C\left(\|b^0\|_{L^2}^2e^{2M_NT^0}+\||\partial_{1}|^\alpha b_n\|_{\dot{H}^{2-2\alpha}}^2+\||\partial_{2}|^\beta b_n\|_{\dot{H}^{2-2\alpha}}^2\right)\|b_n\|_{\dot{H}^{2-2\alpha}}\\
 		\label{3.7}& \leq C\left(\||\partial_{1}|^\alpha b_n\|_{\dot{H}^{2-2\alpha}}^2+\||\partial_{2}|^\beta b_n\|_{\dot{H}^{2-2\alpha}}^2\right)\|b_n\|_{\dot{H}^{2-2\alpha}}+C_N(1+\|b_n\|_{\dot{H}^{2-2\alpha}}^2)
 	\end{align}
 	Collecting the estimates \eqref{3.4}, \eqref{3.5}, \eqref{3.6} and \eqref{3.7}, we obtain
 	\begin{align*}
 		\nonumber	 \frac{d}{dt}\|b_n\|^2_{\dot{H}^{2-2\alpha}}+\||\partial_1|^\alpha b_n\|^2_{\dot{H}^{2-2\alpha}}+\||\partial_2|^\beta b_n\|^2_{\dot{H}^{2-2\alpha}} &\leq C_N(1+\|b_n\|_{\dot{H}^{2-2\alpha}}^2)\\
 		&\hskip0.5cm+C\|b_n\|_{\dot{H}^{2-2\alpha}}\left(\||\partial_{1}|^\alpha b_n\|_{\dot{H}^{2-2\alpha}}^2+\||\partial_{2}|^\beta b_n\|_{\dot{H}^{2-2\alpha}}^2\right).
 	\end{align*}
 	Let $T'_n=\sup\left\{T\in(0,T^0];\ \|b_n\|_{L^\infty_T(\dot{H}^{2-2\alpha})}\leq 2\varepsilon\right\}$.
 	 Since the function $\left(t\mapsto \|b_n\|_{\dot{H}^{2-2\alpha}}\right)$ is continuous in $[0,T^0]$; that  implies
 	the existence of $T'_n$, moreover we have for any $t\in [0,T'_n)$
 	\begin{align*}
 		&\nonumber	 \frac{d}{dt}\|b_n\|^2_{\dot{H}^{2-2\alpha}}+\||\partial_1|^\alpha b_n\|^2_{\dot{H}^{2-2\alpha}}+\||\partial_2|^\beta b_n\|^2_{\dot{H}^{2-2\alpha}} &\leq C_N(1+4\varepsilon^2)
 	+2C\varepsilon\left(\||\partial_{1}|^\alpha b_n\|_{\dot{H}^{2-2\alpha}}^2+\||\partial_{2}|^\beta b_n\|_{\dot{H}^{2-2\alpha}}^2\right).
 	\end{align*}
 	We can choose $\varepsilon<\frac{1}{4C}$, so we get:
 	\begin{align*}
 		2C\varepsilon\left(\||\partial_{1}|^\alpha b_n\|_{\dot{H}^{2-2\alpha}}^2+\||\partial_{2}|^\beta b_n\|_{\dot{H}^{2-2\alpha}}^2\right)\leq \frac{1}{2} \||\partial_{1}|^\alpha b_n\|_{\dot{H}^{2-2\alpha}}^2+ \frac{1}{2} \||\partial_{2}|^\beta b_n\|_{\dot{H}^{2-2\alpha}}^2.
 	\end{align*}
 	Therefore
 	\begin{align*}
 		\nonumber	
 		\frac{d}{dt}\|b_n\|^2_{\dot{H}^{2-2\alpha}}+\frac{1}{2}\||\partial_1|^\alpha b_n\|^2_{\dot{H}^{2-2\alpha}}+\frac{1}{2}\||\partial_2|^\beta b_n\|^2_{\dot{H}^{2-2\alpha}} &\leq C_N(1+4\varepsilon^2).
 	\end{align*}
 	We integrate on $[0, t]$, for any $t\in (0,T'_n)$, we get
 	\begin{align*}
 		\|b_n(t)\|^2_{\dot{H}^{2-2\alpha}}+\int_0^t\||\partial_1|^\alpha b_n(\tau)\|^2_{\dot{H}^{2-2\alpha}}d\tau+\int_{0}^t\||\partial_2|^\beta b_n(\tau)\|^2_{\dot{H}^{2-2\alpha}}d\tau &\leq\|b^0\|^2_{\dot{H}^{2-2\alpha}}+ C_Nt(1+4\varepsilon^2).
 	\end{align*}
 	Taking
 	$$T_0=\frac{\varepsilon^2}{2C_N(1+4\varepsilon^2)}.$$
 	We can show that $T_0<T'_n$, for any $n\in \N$, moreover, for any $t\in [0,T_0]$
 	\begin{equation}
 		\|b_n(t)\|^2_{\dot{H}^{2-2\alpha}}+\int_0^t\||\partial_1|^\alpha b_n(\tau)\|^2_{\dot{H}^{2-2\alpha}}d\tau+\int_{0}^t\||\partial_2|^\beta b_n(\tau)\|^2_{\dot{H}^{2-2\alpha}}d\tau \leq2\varepsilon.
 	\end{equation}
 	So we can conclude that $(b_n)_{n\in \N}$ is uniformly bounded in $C_{T_0}(H^{2-2\alpha}(\R^2))$ and 
 	$$\big(|\partial_{1}|^\alpha b_n\big)_{n},   \big(|\partial_{2}|^\beta b_n\big)_{n}\subset L^2_{T_0}\left(H^{2-2\alpha}(\R^2)\right).$$
 	\\
 	
 	Thanks to these bounds, there remains the passage to the limit of this suite of solutions. This passage to the limit due to the classical argument by combining Ascoli's theorem and Cantor's diagonal process. Finally we get that the system \ref{(AQG)'} admits a solution
 	$b\in \mathcal{X}_{T_0}$.\\
 	
 	Therefore $\theta=a+b$ is solution of system \ref{AQG}  in
 	$ \mathcal{X}_{T_0}.$

 \vskip0.5cm\noindent{\bf $\bullet$ \underline{\large Uniqueness:}} It remains to show that this solution is unique in the space $\mathcal{X}_{T_0}.$ \\
 
 For that let $\theta^1$ and $\theta^2$ be two solutions of \ref{AQG} , $u_{\theta^1}=\mathcal{R}^\perp \theta^1$ and $u_{\theta^2}=\mathcal{R}^\perp \theta^2$. We assume $\omega=\theta^1-\theta^2$, then $\omega$ is solution of the following equation
 	\begin{equation}\label{3.9}
 		\partial_t\omega+|\partial_1|^{2\alpha}\omega+|\partial_2|^{2\beta} \omega+u_1.\nabla\omega+u_\omega.\nabla\theta^2=0,
 	\end{equation}
 	where $u_\omega= \mathcal{R}^\perp \omega=u_{\theta^1}-u_{\theta^2}$.
 	Taken the inner product of \eqref{3.9} with $\omega$, we get
 	\begin{align*}
 		\frac{1}{2}\dfrac{d}{dt}\|\omega(t)\|_{L^2}^2+\||\partial_{1}|^\alpha\omega\|_{L^2}^2+\||\partial_{2}|^\beta\omega\|_{L^2}^2&\leq \left|\left(v.\nabla\theta^2,\omega\right)_{L^2}\right|\\
 		& \leq \left\|v.\nabla\theta^2\right\|_{L^2}\left\|\omega\right\|_{L^2}.
 	\end{align*}
 	Using \cref{Lemma2.3} for $s_1=\alpha<1,\ s_2=1-\alpha<1,$ we get
 	\begin{align*}
 		\left\|v.\nabla\theta^2\right\|_{L^2}
 		&\leq C \||\nabla|^\alpha \omega\|_{L^2} \||\nabla|^\alpha\theta^2\|_{\dot{H}^{2-2\alpha}},
 	\end{align*}
 	and by \cref{Lemma2.7}, we have
 	\begin{align*}
 		\left\|v.\nabla\theta^2\right\|_{L^2}&\leq C \left(\| \omega\|_{L^2}+\||\partial_1|^\alpha \omega\|_{L^2}+\||\partial_2|^\beta \omega\|_{L^2}\right)\||\nabla|^\alpha\theta^2\|_{\dot{H}^{2-2\alpha}}\left\|\omega\right\|_{L^2}\\
 		&\leq \||\partial_1|^\alpha \omega\|_{L^2}+\||\partial_2|^\beta \omega\|_{L^2}+ C \left(1+\||\nabla|^\alpha\theta^2\|_{\dot{H}^{2-2\alpha}}^2\right)\left\|\omega\right\|_{L^2}^2.
 	\end{align*}
 Therefore
 	\begin{align*}
 		\dfrac{d}{dt}\|\omega(t)\|_{L^2}^2\leq C \left(1+\||\nabla|^\alpha\theta^2\|_{\dot{H}^{2-2\alpha}}^2\right)\left\|\omega\right\|_{L^2}^2.
 	\end{align*}
 	Integrate on $[0,t]$, $t\in (0,T_0]$ we get
 	\begin{align*}
 		\|\omega(t)\|_{L^2}^2&\leq C\int_{0}^{t}\left(1+\||\nabla|^\alpha\theta^2(\tau)\|_{\dot{H}^{2-2\alpha}}^2\right)\left\|\omega\right\|_{L^2}^2\ d\tau
 	\end{align*}
 Using Gronwall’s Lemma and the fact $\left(t\mapsto 1+\||\nabla|^\alpha\theta^2(\tau)\|_{\dot{H}^{2-2\alpha}}^2 \right)\in L^1([0,T_0])$ we can deduce that $\omega=0$ in
 $[0, T_0]$ which gives the uniqueness.
 \subsection{\bf Continuity} Finally, we will show that the solution is continuous in $H^{2-2\alpha}(\R^2).$\\
 
 	 First of all, we have $\theta\in C_{T_0}(H^s(\R^2))$, $s<2-2\alpha$. Therefore, if $(s_k)_{k\in \N}$ a positive real sequence such that $$1<s_k<s_{k+1}<2-2\alpha\mbox{ and } \ds s_k\underset{k\rightarrow +\infty}{\longrightarrow}2-2\alpha,$$ then, we have
 	\begin{equation}\label{3.10}
 	\frac{d}{dt}\|\theta(t)\|_{H^{s_k}}^2+2\||\partial_{1}|^\alpha\theta\|_{H^{s_k}}^2+2\||\partial_{2}|^\beta\theta\|_{H^{s_k}}^2+2 \left(|\nabla|^{s_k}\left(u_\theta.\nabla\theta\right),|\nabla|^{s_k}\theta\right)_{L^2}=0.
 	\end{equation}
 	Now, let $0\leq t<t'\leq T_0$, we integer \eqref{3.10} in $[t,t']$, we get
 	\begin{equation}
 		\label{3.11}	\|\theta(t)\|_{H^{s_k}}^2 \leq\|\theta(t')\|_{H^{s_k}}^2+2\int_{t}^{t'}\left(\||\partial_{1}|^\alpha\theta\|_{H^{s_k}}^2d\tau+\||\partial_{2}|^\beta\theta\|_{H^{s_k}}^2\right)d\tau+2\int_{t}^{t'}\left|\left(|\nabla|^{s_k}\left(u_\theta.\nabla\theta\right),|\nabla|^{s_k}\theta\right)_{L^2}\right|d\tau,
 	\end{equation}
 	and 
 	\begin{equation}
 		\label{3.12}	\|\theta(t')\|_{H^{s_k}}^2 \leq\|\theta(t)\|_{H^{s_k}}^2+2\int_{t}^{t'}\left(\||\partial_{1}|^\alpha\theta\|_{H^{s_k}}^2d\tau+\||\partial_{2}|^\beta\theta\|_{H^{s_k}}^2\right)d\tau+2\int_{t}^{t'}\left|\left(|\nabla|^{s_k}\left(u_\theta.\nabla\theta\right),|\nabla|^{s_k}\theta\right)_{L^2}\right|d\tau.
 	\end{equation}
 	
 	Using \cref{Lemma2.5} and \cref{Lemma2.7} and the fact that $u_\theta$ has divergence free and $1<s_k<2-2\alpha$, we get 
 	\begin{align}
 		\nonumber	\int_{t}^{t'}\left|\left(|\nabla|^{s_k}\left(u_\theta.\nabla\theta\right),|\nabla|^{s_k}\theta\right)_{L^2}\right|d\tau&\leq  \int_{t}^{t'}\||\nabla|^{s_k}\left(u_\theta.\nabla\theta\right)-u_\theta.|\nabla|^{s_k}\nabla\theta\|_{L^2}\||\nabla|^{s_k}\theta\|_{L^2}d\tau\\
 		\nonumber	&\leq s_k2^{s_k}C \int_{t}^{t'}\||\nabla|^{s_k+\alpha}\theta\|_{L^2}\||\nabla|^{2-\alpha}\theta\|_{L^2}\|\theta\|_{H^{2-2\alpha}}d\tau\\
 		\nonumber	&\leq C(\alpha) \|\theta\|_{L_{T_0}^\infty(H^{2-2\alpha})} \int_{t}^{t'}\||\nabla|^\alpha\theta\|_{H^{2-2\alpha}}^2d\tau\\
 	\label{3.13}	&\leq C(\alpha)\int_{t}^{t'}\left(\||\partial_{1}|^\alpha\theta\|_{H^{2-2\alpha}}^2+\||\partial_{2}|^\beta\theta\|_{H^{2-2\alpha}}^2\right)d\tau+C(\alpha)(t'-t).
 	\end{align}
 	Then, we collect \eqref{3.11} and \eqref{3.12} with  the inequality \eqref{3.13} to get
\begin{equation}
	\|\theta(t)\|_{H^{s_k}}^2 \leq\|\theta(t')\|_{H^{s_k}}^2+C(\alpha)\int_{t}^{t'}\left(\||\partial_{1}|^\alpha\theta\|_{H^{2-2\alpha}}^2d\tau+\||\partial_{2}|^\beta\theta\|_{H^{2-2\alpha}}^2\right)d\tau+C(\alpha)(t'-t),
\end{equation}
and 
\begin{equation}
\|\theta(t')\|_{H^{s_k}}^2 \leq\|\theta(t)\|_{H^{s_k}}^2+C(\alpha)\int_{t}^{t'}\left(\||\partial_{1}|^\alpha\theta\|_{H^{2-2\alpha}}^2d\tau+\||\partial_{2}|^\beta\theta\|_{H^{2-2\alpha}}^2\right)d\tau+C(\alpha)(t'-t).
\end{equation}
 We have $\ds \lim_{k\rightarrow+\infty}(1+|\xi|^{2s_k})|\widehat{\theta}(\xi,t)|^2= (1+|\xi|^{2(2-2\alpha)})|\widehat{\theta}(\xi,t)|^2$, moreover 
 $$\|\theta(t)\|_{H^{s_k}}^2=\int_{|\xi|\leq 1}(1+|\xi|^{2s_k})|\widehat{\theta}(\xi,t)|^2d\xi+\int_{|\xi|> 1}(1+|\xi|^{2s_k})|\widehat{\theta}(\xi,t)|^2d\xi.$$ 
 The fact that 
 $$\int_{|\xi|\leq 1}(1+|\xi|^{2s_k})|\widehat{\theta}(\xi,t)|^2d\xi\leq 2 \|\theta^0\|_{L^2}^2,$$
 and by Lebesgue's dominated convergence theorem, we get
 	\begin{equation*}
 		\lim_{k\rightarrow+\infty}\int_{|\xi|\leq 1}(1+|\xi|^{2s_k})|\widehat{\theta}(\xi,t)|^2d\xi=\int_{|\xi|\leq 1}(1+|\xi|^{2(2-2\alpha)})|\widehat{\theta}(\xi,t)|^2d\xi.
 	\end{equation*}
 Moreover, we have 
 $$ (1+|\xi|^{2s_k})|\widehat{\theta}(\xi,t)|^2\leq (1+|\xi|^{2s_{k+1}})|\widehat{\theta}(\xi,t)|^2,\quad\forall |\xi|>1,$$
 therefore, by the Monotonic Convergence Theorem, we obtain
  	\begin{equation}
  	\lim_{k\rightarrow+\infty}\int_{|\xi|> 1}(1+|\xi|^{2s_k})|\widehat{\theta}(\xi,t)|^2d\xi=\int_{|\xi|> 1}(1+|\xi|^{2(2-2\alpha)})|\widehat{\theta}(\xi,t)|^2d\xi.
  \end{equation}
Finally, we have for any $t\in [0,T_0]$, $\ds \lim_{k\rightarrow+\infty}\|\theta(t)\|_{H^{s_k}}=\|\theta(t)\|_{H^{2-2\alpha}}$, which implies
\begin{equation}
	\|\theta(t)\|_{H^{2-2\alpha}}^2 \leq\|\theta(t')\|_{H^{2-2\alpha}}^2+C(\alpha)\int_{t}^{t'}\left(\||\partial_{1}|^\alpha\theta\|_{H^{2-2\alpha}}^2d\tau+\||\partial_{2}|^\beta\theta\|_{H^{2-2\alpha}}^2\right)d\tau+C(\alpha)(t'-t),
\end{equation}
and 
\begin{equation}
	\|\theta(t')\|_{H^{2-2\alpha}}^2 \leq\|\theta(t)\|_{H^{2-2\alpha}}^2+C(\alpha)\int_{t}^{t'}\left(\||\partial_{1}|^\alpha\theta\|_{H^{2-2\alpha}}^2d\tau+\||\partial_{2}|^\beta\theta\|_{H^{2-2\alpha}}^2\right)d\tau+C(\alpha)(t'-t).
\end{equation}
After that, we pass $\ds\limsup_{t\rightarrow t'}$ and $\ds\limsup_{t'\rightarrow t}$ for the last equations, then, we get
 	\begin{equation}
 		\ds\limsup_{t\rightarrow t'}	\|\theta(t)\|_{H^{2-2\alpha}}\leq 	\|\theta(t')\|_{H^{2-2\alpha}}\quad\mbox{ and } \quad \ds\limsup_{t'\rightarrow t}	\|\theta(t')\|_{H^{2-2\alpha}}\leq 	\|\theta(t)\|_{H^{2-2\alpha}},
 	\end{equation}
 	moreover $\theta(t)\rightharpoonup \theta(t')$ if $t\rightarrow t'$ and $\theta(t)\rightharpoonup \theta(t')$ if $t'\rightarrow t$. Therefore by \cref{Lemma2.1} we have $$\lim\limits_{t\rightarrow t'}\|\theta(t)-\theta(t')\|_{H^{2-2\alpha}}=0\quad\mbox{ and }\quad\lim\limits_{t'\rightarrow t}\|\theta(t')-\theta(t)\|_{H^{2-2\alpha}}=0,$$ which implies that $\left(t\mapsto\theta(t)\right)$ is continue in $[0,T_0]$ in $H^{2-2\alpha}(\R^2)$.\\
 	
 	So we can conclude that $$\theta\in C([0,T_0],H^{2-2\alpha}(\R^2)).$$ 

 \subsection{\bf Global solution} In this section, we prove if $\|\theta^0\|_{\dot{H}^{2-2\alpha}}<c$, we get a global solution in $C(\R^+, H^{2-2\alpha}(\R^2))$ satisfying \eqref{1.2}. \\
 
 	Let $\theta\in C([0,T^\ast), H^{2-2\alpha}(\R^2))$ be a maximal solution of the system \ref{AQG}, using \cref{Lemma2.5} and  the fact $\dive u_\theta=0$, we infer that
 	\begin{align}
 		\nonumber	 \frac{1}{2}\frac{d}{dt}\|\theta(t)\|^2_{H^{2-2\alpha}}+\||\partial_{1}|^\alpha\theta\|_{H^{2-2\alpha}}^2+\||\partial_{2}|^\beta\theta\|_{H^{2-2\alpha}}^2&\leq \left|\left(|\nabla|^{2-2\alpha}(u_\theta.\nabla\theta),|\nabla|^{2-2\alpha}\theta\right)_{L^2}\right|\\
 		\label{3.20}	&\leq C \||\nabla|^{2-\alpha}\theta\|_{L^2}^2\||\nabla|^{2-2\alpha}\theta\|_{L^2}
 	\end{align}
 	By interpolation, we get 
 	\begin{equation*}
 		\||\nabla|^{2-\alpha}\theta\|_{L^2}^2\leq C\left( \||\partial_1|^{\alpha}\theta\|_{H^{2-2\alpha}}^2+\||\partial_2|^{\beta}\theta\|_{H^{2-2\alpha}}^2\right),
 	\end{equation*}
 which implies
 	\begin{equation}
 		\frac{d}{dt}\|\theta(t)\|^2_{H^{2-2\alpha}}+2\||\partial_{1}|^\alpha\theta\|_{H^{2-2\alpha}}^2+2\||\partial_{2}|^\beta\theta\|_{H^{2-2\alpha}}^2\leq C \left(\||\partial_1|^{\alpha}\theta\|_{H^{2-2\alpha}}^2+ \||\partial_2|^{\beta}\theta\|_{H^{2-2\alpha}}^2\right)\|\theta\|_{\dot{H}^{2-2\alpha}}
 	\end{equation}
 	Let $T_1=\sup\left\{T\in(0,T^\ast);\ \|\theta\|_{L^\infty_T(\dot{H}^{2-2\alpha})}\leq 2\|\theta^0\|_{\dot{H}^{2-2\alpha}}\right\}$, then for any $t\in [0,T_1)$, we have
 	\begin{equation}
 		\frac{d}{dt}\|\theta(t)\|^2_{H^{2-2\alpha}}+2\||\partial_{1}|^\alpha\theta\|_{H^{2-2\alpha}}^2+2\||\partial_{2}|^\beta\theta\|_{H^{2-2\alpha}}^2\leq 2C \left(\||\partial_1|^{\alpha}\theta\|_{H^{2-2\alpha}}^2+ \||\partial_2|^{\beta}\theta\|_{H^{2-2\alpha}}^2\right)\|\theta^0\|_{\dot{H}^{2-2\alpha}}.
 	\end{equation}
 	Taking  $c=\dfrac{1}{2 C}$, then, if we have  $\|\theta^0\|_{H^{2-2\alpha}}<c$, so for any $t\in  [0,T_1)$
 	\begin{equation}
 		\frac{d}{dt}\|\theta(t)\|^2_{H^{2-2\alpha}}+\||\partial_{1}|^\alpha\theta\|_{H^{2-2\alpha}}^2+\||\partial_{1}|^\alpha\theta\|_{H^{2-2\alpha}}^2\leq 0,
 	\end{equation}
 	which implies
 	\begin{equation*}
 		\|\theta(t)\|^2_{H^{2-2\alpha}}+\int_{0}^{t}\||\partial_{1}|^\alpha\theta(\tau)\|_{H^{2-2\alpha}}^2d\tau+\int_{0}^{t}\||\partial_{1}|^\alpha\theta(\tau)\|_{H^{2-2\alpha}}^2d\tau\leq \|\theta^0\|^2_{H^{2-2\alpha}}.
 	\end{equation*}
Finally, we get $T = T^\ast$. Hence $T^\ast=+\infty$, which complete the proof.
\section{\bf Proof of \cref{theorem1.3}}
 	The proof is done in two steps. In the first, we prove
 	\begin{equation}
 		\lim\limits_{t\rightarrow+\infty}\|\theta(t)\|_{\dot{H}^{2-2\alpha}}=0,
 	\end{equation}
 	where we use
 	\begin{equation*}
 		\|\theta(t)\|^2_{H^{2-2\alpha}}+\int_{0}^{t}\||\partial_{1}|^\alpha\theta(\tau)\|_{H^{2-2\alpha}}^2d\tau+\int_{0}^{t}\||\partial_{1}|^\alpha\theta(\tau)\|_{H^{2-2\alpha}}^2d\tau\leq \|\theta^0\|^2_{H^{2-2\alpha}},\quad \forall t\geq 0.
 	\end{equation*}
 	In the second step, we prove
 	\begin{equation}
 		\lim\limits_{t\rightarrow+\infty}\|\theta(t)\|_{L^2}=0.
 	\end{equation}
 	\textbf{\underline{First step:}} By interpolation theorem we have
	$$\|\theta(t)\|^2_{\dot{H}^{2-2\alpha}}\leq C\left(\||\partial_{1}|^{\alpha}\theta(t)\|^2_{H^{2-2\alpha}}+\||\partial_{2}|^{\beta}\theta(t)\|_{H^{2-2\alpha}}^2\right).$$
 	We integer over $\R^+$ we get
 	\begin{equation}
 		\int_{0}^{+\infty}\|\theta(t)\|^2_{\dot{H}^{2-2\alpha}}dt\leq C\int_{0}^{+\infty}\||\partial_{1}|^{\alpha}\theta(t)\|_{H^{2-2\alpha}}^2dt+C\int_{0}^{+\infty}\||\partial_{2}|^{\beta}\theta(t)\|_{H^{2-2\alpha}}^2dt\leq C\|\theta^0\|_{H^{2-2\alpha}}^2.
 	\end{equation}
 	Let $\varepsilon>0$ and $E(\varepsilon)=\{t\geq 0;\ \|\theta(t)\|_{\dot{H}^{2-2\alpha}}>\varepsilon\}$, then
 	\begin{equation}
 		\varepsilon^2\lambda(E)=\int_{E}\|\theta(t)\|^2_{\dot{H}^{2-2\alpha}}dt\leq \int_{0}^{+\infty}\|\theta(t)\|^2_{\dot{H}^{2-2\alpha}}dt\leq C\|\theta^0\|_{H^{2-2\alpha}}^2.
 	\end{equation}
 	where $\lambda(E)$ is the Lebesgue measure of $E$, then $\lambda(E)<+\infty$. So for $r>0$, there exists $t_0\in [0,\lambda(E)+r]$ such that $t_0\notin E$ and
 	$$\|\theta(t_0)\|_{\dot{H}^{2-2\alpha}}\leq\varepsilon.$$
 	Thus, $\ds\lim_{t\rightarrow +\infty}\|\theta(t)\|_{\dot{H}^{2-2\alpha}}=0$.\\
 	\textbf{\underline{Second step:}} Let $\delta$ a strictly positive real number strictly less than 1, we define the operators
 	$A_\delta(D)$ and $B_\delta(D)$, respectively, by the following:
 	\begin{equation}
 		\begin{array}{l}
 			A_\delta(D)f=\mathcal{F}^{-1}\left(\xi\mapsto\chi_{\mathcal{P}(0,\delta)}(\xi)\mathcal{F}(f)(\xi)\right),\\
 			B_\delta(D)f=\mathcal{F}^{-1}\left(\xi\mapsto(1-\chi_{\mathcal{P}(0,\delta)})(\xi)\mathcal{F}(f)(\xi)\right),
 		\end{array}
 	\end{equation}
 where
 $$\mathcal{P}(0,\delta)=\left\{\xi\in\R^2:\ A(\xi)\leq \delta\right\}.$$
 	We define $w_\delta=A_\delta(D)\theta$ and $v_\delta=B_\delta(D)\theta$; $\mathcal{F}(\theta)=\mathcal{F}(w_\delta)+\mathcal{F}(v_\delta)$. Then,
 	\begin{equation}\label{4.6}
 		\partial_tw_\delta+|\partial_1|^{2\alpha}w_\delta+|\partial_2|^{2\beta}w_\delta+A_\delta(D)( u_\theta.\nabla\theta)=0
 	\end{equation}
 	and
 	\begin{equation}
 		\partial_tv_\delta+|\partial_1|^{2\alpha}v_\delta+|\partial_2|^{2\beta}v_\delta+B_\delta(D)( u_\theta.\nabla\theta)=0
 	\end{equation}
 	Taking the scalar product of \eqref{4.6} equation with $w_\delta$, we get
 	\begin{align*}
 		\frac{1}{2}	\frac{d}{dt}\|w_\delta(t)\|_{L^2}^2&\leq \left|\left(A_\delta(D)(u_\theta.\nabla\theta),w_\delta\right)_{L^2}\right|\\
 		&\leq \int_{A(\xi)\leq \delta}|\xi||\widehat{u_\theta\otimes \theta}(\xi)||\widehat{w_\delta}(\xi)|d\xi\\
 		&\leq \int_{A(\xi)\leq \delta}|\xi|^{2-2\beta}|\xi|^{2\beta-1}|\widehat{u_\theta\otimes \theta}(\xi)||\widehat{w_\delta}(\xi)|d\xi\\
 		&\leq C(\alpha,\beta) \int_{A(\xi)\leq \delta}\left(A(\xi)^{\frac{1}{2\alpha}}+A(\xi)^{\frac{1}{2\beta}}\right)^{2-2\beta}|\xi|^{2\beta-1}||\widehat{u_\theta\otimes \theta}(\xi)||\widehat{w_\delta}(\xi)|d\xi\\
 		&\leq C(\alpha,\beta)  \delta^{1-\beta}\int_{\R^2}|\xi|^{2\beta-1}|\widehat{u_\theta\otimes \theta}(\xi)||\widehat{w_\delta}(\xi)|d\xi\\
 		&\leq C(\alpha,\beta) \delta^{1-\beta}\|u_\theta\otimes \theta\|_{\dot{H}^{2\beta-1}}\|w_\delta\|_{L^2}\\
 		&\leq C(\alpha,\beta)  \delta^{1-\beta}\||\nabla|^{\beta}\theta\|_{L^2}^2\|w_\delta\|_{L^2}.
 	\end{align*}
 	By interpolation Theorem, $\alpha\leq\beta\leq 2-\alpha$ we have
 	\begin{align*}
 		\||\nabla|^{\beta}\theta\|_{L^2}&\leq \||\partial_{1}|^\beta\theta\|_{L^2}+\||\partial_{2}|^\beta\theta\|_{L^2}\\
 		&\leq \||\partial_{1}|^\alpha\theta\|_{L^2}+ \||\partial_{1}|^\alpha\theta\|_{\dot{H}^{2-2\alpha}}+\||\partial_{2}|^\beta\theta\|_{L^2}
 	\end{align*}
 	Using $\|w_\delta\|_{L^2}\leq \|\theta\|_{L^2}\leq \|\theta^0\|_{L^2}$, and integer in $[0,t]$ we get
 	\begin{align*}
 		\|w_\delta(t)\|_{L^2}^2
 		&\leq 	\|w_\delta^0\|_{L^2}^2+C(\alpha,\beta)\delta^{1-\beta}\|\theta^0\|_{L^2}\int_{0}^{t}\left(\||\partial_{1}|^\alpha\theta\|_{H^{2-2\alpha}}^2+\||\partial_{2}|^\beta\theta\|_{L^2}^2\right)d\tau\\
 		&\leq \|w_\delta^0\|_{L^2}^2+ C(\alpha,\beta)\delta^{1-\beta}\|\theta^0\|_{H^{2-2\alpha}}^3\underset{\delta\rightarrow0^+}{\longrightarrow}0.
 	\end{align*}
 	Therefore
 	\begin{equation*}
 		\lim\limits_{\delta\rightarrow 0^+}\|w_\delta(t)\|_{L^\infty(\R^+,L^2)}=0.
 	\end{equation*}
So for any $\varepsilon>0$, there exists $\delta_0>0$ such that 
 	\begin{equation}\label{4.8}
 		\|w_{\delta_0}(t)\|_{L^\infty(\R^+,L^2(\R^2))}<\dfrac{\varepsilon}{2}.
 	\end{equation}
 	On the other hand, we deduce that
 	\begin{align*}
 		\|v_\delta\|^2_{L^2}&\leq \int_{A(\xi)>\delta}\frac{A(\xi)}{A(\xi)}|\widehat{\theta}(\xi)|^2d\xi\\
 		& \leq \frac{1}{\delta}\int_{\R^2}\left(|\xi_1|^{2\alpha}+|\xi_2|^{2\beta}\right)|\widehat{\theta}(\xi)|^2d\xi\\
 		&\leq \frac{1}{\delta}\left(\||\partial_{1}|^\alpha\theta\|^2+\||\partial_{2}|^\beta\theta\|^2\right),
 	\end{align*}
 which implies 
 \begin{equation*}
 	\int_{0}^{+\infty}\|v_\delta\|^2_{L^2}d\tau\leq \frac{1}{{\delta}}\int_{0}^{+\infty}\left(\||\partial_{1}|^\alpha\theta\|^2+\||\partial_{2}|^\beta\theta\|^2\right)d\tau\leq \frac{1}{{\delta}}\|\theta^0\|^2.
 \end{equation*}
 	Let $S_\varepsilon(\delta_0)=\{t\geq 0; \|v_{\delta_0}\|_{L^2}>\frac{\varepsilon}{2}\}$, then
 	\begin{align*}
 		\left(\frac{\varepsilon}{2}\right)^2\lambda(S_\varepsilon(\delta_0))&\leq \int_{S_\varepsilon(\delta_0)}\|v_{\delta_0}\|_{L^2}^2dt\\
 		&\leq \int_{0}^{+\infty}\|v_{\delta_0}\|_{L^2}^2dt\\
 		&\leq \frac{1}{{\delta_0}}\|\theta^0\|_{L^2}^2,
 	\end{align*}
 	We pose
 	$$T_{\varepsilon}=\left(\frac{2}{\varepsilon}\right)^2\frac{1}{{\delta_0}}\|\theta^0\|_{L^2}^2<+\infty,$$
 	then $\lambda(S_\varepsilon(\delta_0))\leq T_{\varepsilon}$. So there exists $t_2\in [0,T_\varepsilon+1]\setminus S_\varepsilon(\delta_0)$ such
 	\begin{equation}\label{4.9}
 		\|v_{\delta_0}(t_2)\|_{L^2}\leq \frac{\varepsilon}{2}.
 	\end{equation}
 	By the equation \eqref{4.8} and \eqref{4.9}, we get
 	\begin{equation*}
 		\|\theta(t_2)\|_{L^2}\leq \varepsilon.
 	\end{equation*}
 	Thus, $\lim\limits_{t\rightarrow+\infty}\|\theta(t)\|_{L^2}=0$, and this finishes the proof in this case.\hfill$\blacksquare$
	

\medskip

\begin{thebibliography}{9}
		\bibitem{HB} Brezis, H. (1983). Analyse fonctionnelle. Théorie et applications.
	
	\bibitem{BH} 
	Bahouri, H., Chemin, J. Y., \& Danchin, R. (2011). Fourier analysis and nonlinear partial differential equations (Vol. 343, pp. 523-pages). Berlin: Springer.
	
	\bibitem{JP} Pedlosky, J. (1987). Geophysical fluid dynamics (Vol. 710). New York: springer.
	
		\bibitem{DC}  \href{https://www.jstor.org/stable/121037?casa_token=R9hHmwlVh9kAAAAA%3AH9oUnpoRsjzfnWxHN3Ps2NSBoDLr8EWvj0ON4FsTkG_NkZdGWbHjuvyB6tQfFGYCT0WMQzTcnanPgNtFynGUbHPS8LQeFpUeuPxu8gkT7V8XB13XWT3t}{Cordoba, D. (1998). Nonexistence of simple hyperbolic blow-up for the quasi-geostrophic equation. Annals of Mathematics, 148(3), 1135-1152.}
	
	
	\bibitem{CP} \href{https://onlinelibrary.wiley.com/doi/abs/10.1002/cpa.3160380605}{Constantin, P., Lax, P. D., \& Majda, A. (1985). A simple one‐dimensional model for the three‐dimensional vorticity equation. Communications on pure and applied mathematics, 38(6), 715-724.}
	
	\bibitem{CP1} 
	\href{https://iopscience.iop.org/article/10.1088/0951-7715/7/6/001/meta?casa_token=I5R_JR3my10AAAAA:AutKyYCQDRLhAOXPh_y0Rdm1abaY5XtNZ-bLKQpCfWP57GomC-7OREv5wH-fSSvxD6sWnwIoElJbXQ}{Constantin, P., Majda, A. J., \& Tabak, E. (1994). Formation of strong fronts in the 2-D quasigeostrophic thermal active scalar. Nonlinearity, 7(6), 1495.}

\bibitem{CW} \href{https://epubs.siam.org/doi/abs/10.1137/S0036141098337333}{Constantin, P., \& Wu, J. (1999). Behavior of solutions of 2D quasi-geostrophic equations. SIAM journal on mathematical analysis, 30(5), 937-948.}

\bibitem{CV} \href{https://link.springer.com/article/10.1007/s00039-012-0172-9}{Constantin, P., \& Vicol, V. (2012). Nonlinear maximum principles for dissipative linear nonlocal operators and applications. Geometric And Functional Analysis, 22(5), 1289-1321.}

\bibitem{CZ} \href{https://www.sciencedirect.com/science/article/pii/S0362546X06004731?casa_token=yspRhX-47Y4AAAAA:GjtEZg65n2huqrWKkKceMJGMCYpOve36_ZUxoqUAYFOFqKjJt_CifxyaHofeVW9UH-dLeC3Kpwq3}{Chen, Q., \& Zhang, Z. (2007). Global well-posedness of the 2D critical dissipative quasi-geostrophic equation in the Triebel–Lizorkin spaces. Nonlinear Analysis: Theory, Methods \& Applications, 67(6), 1715-1725.}

	\bibitem{KNV} \href{https://link.springer.com/article/10.1007/s00222-006-0020-3}{Kiselev, A., Nazarov, F., \& Volberg, A. (2007). Global well-posedness for the critical 2D dissipative quasi-geostrophic equation. Inventiones mathematicae, 167(3), 445-453.}

\bibitem{ZF} \href{https://link.springer.com/article/10.1007/s11425-007-0002-y}{Zhang, Z. F. (2007). Global well-posedness for the 2D critical dissipative quasi-geostrophic equation. Science in China Series A: Mathematics, 50(4), 485-494.}
	
\bibitem{MH} \href{https://link.springer.com/article/10.1007/s00220-006-0023-3}{Miura, H. (2006). Dissipative quasi-geostrophic equation for large initial data in the critical Sobolev space. Communications in mathematical physics, 267(1), 141-157.}
	
\bibitem{JN} 
\href{https://www.jstor.org/stable/24902684?casa_token=DIeQk7C9CaQAAAAA%3A9ycJn-soIX7uX0f_4ux4-Oy9eCBUxzBwERhPNDZUt90nBr1UPXiK1RFMdrN70bhV3Pt-VGzb764x0y3n68IANc37oGKr_rAR0FA4FIOsrLqmyw-nyy9F&seq=1}{Ju, N. (2007). Dissipative 2D quasi-geostrophic equation: local well-posedness, global regularity and similarity solutions. Indiana University mathematics journal, 187-206.}


\bibitem{YZ} \href{https://iopscience.iop.org/article/10.1088/1361-6544/ab41e6/meta}{Ye, Z. (2019). On the global regularity for the anisotropic dissipative surface quasi-geostrophic equation. Nonlinearity, 33(1), 72.}

	
	


\bibitem{YZ1} \href{https://arxiv.org/abs/2012.02956}{Ye, Z. (2020). Global regularity and time decay for the SQG equation with anisotropic fractional dissipation. arXiv preprint arXiv:2012.02956.}


\end{thebibliography}
\end{document}